\documentclass{amsart}
\usepackage{amssymb}
\title[Bimonial Character Sums]{Evaluating Binomial Character Sums
Modulo Powers of two}

\author{Vincent Pigno, Chris  Pinner, \and Joe Sheppard}
\address{ Department of Mathematics\\
         Kansas State University\\
         Manhattan, KS 66506}
\email{pignov@math.ksu.edu, pinner@math.ksu.edu \& jnsheppa@ksu.edu }
\thanks{The first and third  authors acknowledge support of K-State's I-Center and Arts \& Sciences Undergraduate Research Scholarships respectively.}

\date{\today}

\newcommand{\be}{\begin{equation}}
\newcommand{\ee}{\end{equation}}
\newcommand{\ba}{\begin{align}}
\newcommand{\ea}{\end{align}}

\begin{document}

\begin{abstract} We show that for any mod $2^m$ characters, $\chi_1, \chi_2,$   the  complete exponential  sum,
$$
\sum_{x=1}^{2^m}\chi_1(x) \chi_2(Ax^k+B),
$$
has a simple explicit evaluation.
\end{abstract}

\maketitle

\newtheorem{theorem}{Theorem}[section]
\newtheorem{corollary}{Corollary}[section]
\newtheorem{lemma}{Lemma}[section]
\newtheorem{conjecture}{Conjecture}[section]

\section{Introduction}
Suppose that  $\chi_1$ and $\chi_2$  are mod $2^m$ multiplicative characters with $\chi_2$ primitive mod $2^m$, $m\geq 3$.  We  are interested here in evaluating  the complete character sum
$$ S=\sum_{x=1}^{2^m}\chi_1(x)\chi_2(Ax^k+B). $$
Writing $\chi (Ax^{K} +Bx^{L})=\chi^{L}(x) \chi (Ax^{K-L}+B)$ these sums of course include the binomial character sums.
Cases where one can explicitly evaluate an exponential or character sum are unusual and therefore worth investigating.  

In \cite{PPP} we considered the 
corresponding result for mod $p^m$ characters with $p\geq 3$ and $m$ sufficiently large, using reduction techniques of Cochrane \& Zheng \cite{c}  (see also \cite{cz} \& \cite{cz2}).  Though their results are stated for odd primes
the approach  can often be adapted  for $p=2$ as we showed for twisted monomial exponential sums in \cite{PP}.
When $k=1$  and $A=-1$, $B=1$ the mod $p^m$ sum is  the classical Jacobi sum (though uninterestingly  zero if $p=2$). See \cite{BerndtBk} or \cite{LidlNiedBk} for an extensive treatment of mod $p$ Jacobi sums  and their generalizations over $\mathbb F_{p^m}$. 
In \cite{Long} we treated the $k=1$ case for general $A$, $B$,  including when  $p=2$,  along with generalizations of the multivariable Jacobi sums considered in \cite{Wenpeng2}.

Plainly  $S=0$ if $A$ and $B$ 
are not of opposite parity (otherwise $x$ or $Ax^k+B$ will be even and the individual terms will all be zero).  We assume here that $A$ is even and  $B$ is odd  and  write
$$ A=2^nA_1,\; n>0,   \hspace{2ex} k=2^tk_1, \hspace{2ex} 2\nmid A_1k_1B. $$
If $B$ is even and $A$ odd we can use  $x\mapsto x^{-1}$  to write $S$ in the form
$$ S= \sum_{x=1}^{2^m}  \overline{\chi_1}\overline{\chi}_2^k (x) \chi_2 (Bx^k+A). $$
Since $\mathbb Z_{2^m}^*=<-1,5>$,  the characters $\chi_1$, $\chi_2$ are completely determined by their values on $-1$ 
and 5. Since $5$ has order $2^{m-2}$ mod $2^m$ we  can define  integers $c_1$, $c_2$  with
$$ \chi_i (5) =e_{2^{m-2}}(c_i),\hspace{2ex} 1\leq c_i \leq 2^{m-2},$$
where $e_n(x):=e^{2\pi  i x /n}$. Since $\chi_2$ is primitive we have $2\nmid c_2$.  We define  the odd integers  $R_i$, $i\geq 2$, by
\be  \label{defR}   5^{2^{i-2}}= 1 +R_i2^i. \ee
Defining 
$$ N:= \begin{cases}    \lceil \frac{1}{2}(m-n)\rceil, & \text{ if $m-n>2t+4$}, \\
 t+2, & \text{if $t+2\leq m-n \leq 2t+4$,} \end{cases} $$
and
\be \label{char0}      C(x):=c_1(Ax^k+B) + c_2 Akx^kR_{N}R^{-1}_{N+n} \ee
(here and throughout the paper $y^{-1}$ denotes the inverse of $y$ mod $2^m$) 
it transpires that the  sum $S$ will be zero unless there is a solution $x_0$ to the
characteristic equation
\be  \label{defchar}  
  C(x_0) \equiv 0 \text{ mod } 2^{\lfloor \frac{1}{2}(m+n)\rfloor +t}, \ee
with $2\nmid x_0(Ax_0^k+B)$, when $m-n>2t+4$, and a solution to $C(1)$ or $C(-1)\equiv 0$ mod $2^{m-2}$ when $t+2\leq m-n\leq 2t+4$.

\begin{theorem}\label{maintheorem}

Suppose that $m-n \geq   t+2$. The sum $S=0$ unless $c_1=2^{n+t}c_3,$  with  $2\nmid c_3$, and  $\chi_1(-1)=1$ when $k$ is even,  and the characteristic equation \eqref{defchar} has an odd  solution $x_0$  when $m-n > 2t+4$. Assuming these  conditions do hold.

When $m-n>2t+4$,
$$ S=2^{\frac{1}{2}(m+n)+t+\min\{1,t\} } \chi_1(x_0)\chi_2(Ax_0^k+B) \begin{cases}  1 , & \text{ if $m-n$ is even, } \\
\omega^{h}\left(\frac{2}{h}\right), &  \text{ if $m-n$ is odd, }
\end{cases} $$
where $\left(\frac{2}{x}\right)$ is the Jacobi symbol,  $\omega =e^{\pi i/4}$, $C(x_0)=\lambda 2^{\lfloor \frac{1}{2}(m+n)\rfloor +t}$  for some integer $\lambda$ and $h:= 2\lambda +(k_1-1) +(2^n-1)c_3$.

When  $t+3<  m-n\leq 2t+4$, 
$$ S = \begin{cases} 2^{m-1} \chi_2(A+B), & \text{ if $k$ is even and $C(1)\equiv 0$ mod $2^{m-2}$, } \\ 2^{m-2}\chi_2(A+B),  & \text{ if  $k$ is odd  and $C(1)\equiv 0$ mod $2^{m-2}$, } \\ 2^{m-2}\chi_1(-1)\chi_2(-A+B),  & \text{ if  $k$ is odd and $C(-1)\equiv 0$ mod $2^{m-2}$, } \\
 0, & \text{ otherwise.} \end{cases} $$

When $m-n=t+3$,
$$ S=\begin{cases}  2^{m-1} \chi_2(A+B), & \text{ if $k$ is even and $\chi_1(5)=\pm 1$, $\chi_1(-1)=1$,}  \\
 2^{m-2} \left(\chi_2(A+B) +\chi_1(-1) \chi_2 (-A+B)\right),  & \text{ if $k$ is odd and $\chi_1(5)=\pm 1$,} \\
 0, & \text{ otherwise.} \end{cases} $$

When $m-n=t+2$,
$$ S=\begin{cases}  2^{m-1} \chi_2(A+B), & \text{ if $k$ is even and $\chi_1=\chi_0$ or  $k$ is odd and $\chi_1=\chi_4$,}  \\
 0, & \text{ otherwise,} \end{cases} $$
where $\chi_0$ is the principal character mod $2^m$ and $\chi_4$  is the mod $2^m$ character induced by the non-trivial character mod 4
(i.e. $\chi_4(x)=\pm 1$ as $x\equiv \pm 1$ mod 4 respectively).

\end{theorem}
Note that the restriction $m-n \geq  t+2$ is quite  natural; for  $m-n < t+2$ 
the odd $x$ have $Ax^k+B\equiv A+B$ mod $2^m$ and $S=\chi_2(A+B)\sum_{x=1}^{2^m}\chi_1(x)=2^{m-1}\chi_2(A+B)$
if $\chi_1=\chi_0$  and zero otherwise.

Our original assumption that $\chi_2$ is primitive is also reasonable; if $\chi_1$ and $\chi_2$
are both imprimitive then one should  reduce the modulus, while if $\chi_1$ is primitive 
and $\chi_2$ imprimitive then $S=0$ 
(if $\chi_1$ is primitive then  $u=1+2^{m-1}$ must have $\chi_1(u)=-1$, since  $x+2^{m-1}\equiv ux$ mod $2^m$ for any odd $x$,  and $x\mapsto xu$ gives $S=\chi_1(u)S$
when $\chi_2$ is imprimitive).

\section{Proof}

\noindent
{\bf Initial  decomposition}

Observing that $\pm 5^{\gamma}$,  $\gamma =1,\ldots, 2^{m-2}$, gives a reduced residue system mod $2^m$
and writing
$$ S(A) :=\sum_{\gamma=1}^{2^{m-2}} \chi_1(5^{\gamma}) \chi_2(A5^{\gamma k}+B), $$
if $k$ is even we have 
\be \label{trans1}  S= (1+\chi_1 (-1)) S(A)=\begin{cases}  0, & \text{ if $\chi_1(-1)=-1$,} \\ 2S(A), & \text{ if $\chi_1(-1)=1$,} \end{cases} \ee
and if $k$ is odd
\be \label{trans2} S=S(A)+\chi_1(-1) S(-A). \ee

\vspace{2ex}

\noindent
{\bf  Large $m$ values: $m > n+2t+4$}

If  $I_1$ is an interval of length $ 2^{\lceil \frac{m-n}{2}\rceil -t -2}$ then plainly
$$ \gamma = u 2^{\lceil \frac{m-n}{2}\rceil -t -2} +v, \;\; v\in I_1,\;\; u\in I_2:=\left[1, 2^{\lfloor \frac{m+n}{2}\rfloor +t}\right], $$
runs through a complete set of residues mod $2^{m-2}$.   Hence, writing $h(x):=Ax^k+B$ and noting that $2\nmid h(5^v)$, 
\begin{align*} S(A) & =\sum_{v\in I_1}  \chi_1(5^v) \sum_{u\in I_2} \chi_1\left(5^{ u 2^{\lceil \frac{m-n}{2}\rceil -t -2}}\right) \chi_2 \left( A5^{vk} 5^{k    u 2^{\lceil \frac{m-n}{2}\rceil -t -2}}+B\right)  \\  & = \sum_{v\in I_1}  \chi_1(5^v) \chi_2(h(5^v))  \sum_{u\in I_2} \chi_1\left(5^{ u 2^{\lceil \frac{m-n}{2}\rceil -t -2}}\right)\chi_2(W)
\end{align*}
where 
\begin{align*}  W &  =h(5^v)^{-1} A5^{vk} \left(  5^{k    u 2^{\lceil \frac{m-n}{2}\rceil -t -2}}-1\right) +1 .\\
\end{align*}
Since $n+ 2\lceil \frac{m-n}{2}\rceil \geq m$ and $2\lceil \frac{m+n}{2}\rceil \geq m$ we have
\begin{align*}
W & = A_15^{vk}h(5^v)^{-1} 2^n \left( \left(1+R_{\lceil \frac{m-n}{2}\rceil }2^{\lceil \frac{m-n}{2}\rceil }\right)^{uk_1}-1\right)  +1\\
& \equiv 1 + A_15^{vk}h(5^v)^{-1}uk_1 R_{\lceil \frac{m-n}{2}\rceil }2^{\lceil \frac{m+n}{2}\rceil} \text{ mod }2^m \\
 & \equiv \left( 1 + R_{\lceil \frac{m+n}{2}\rceil }2^{\lceil \frac{m+n}{2}\rceil }\right)^{A_15^{vk}h(5^v)^{-1}uk_1 R_{\lceil \frac{m-n}{2}\rceil }R_{\lceil \frac{m+n}{2}\rceil }^{-1}}  \text{ mod } 2^m \\
 & = 5^{ A_15^{vk}h(5^v)^{-1}uk_1 R_{\lceil \frac{m-n}{2}\rceil }R_{\lceil \frac{m+n}{2}\rceil }^{-1}2^{\lceil \frac{m+n}{2}\rceil -2}} \\
 & =  5^{ A5^{vk}h(5^v)^{-1}uk R_{\lceil \frac{m-n}{2}\rceil }R_{\lceil \frac{m+n}{2}\rceil }^{-1}2^{\lceil \frac{m-n}{2}\rceil -t-2}}   .
\end{align*}
So we can write 
$$
\sum_{u\in I_2} \chi_1\left(5^{ u 2^{\lceil \frac{m-n}{2}\rceil -t -2}}\right)\chi_2(W)  = \sum_{u\in I_2}  e_{2^{\lfloor \frac{m+n}{2}\rfloor +t}}\left(u \left(c_1+ c_2 A5^{vk}h(5^v)^{-1}k R_{\lceil \frac{m-n}{2}\rceil }R_{\lceil \frac{m+n}{2}\rceil }^{-1}\right)\right) ,$$
which equals $ 2^{\lfloor \frac{m+n}{2}\rfloor +t}$  for the $v$ with
\be \label{char2}  c_1 h(5^v) + c_2 A5^{vk} k R_{\lceil \frac{m-n}{2}\rceil }R_{\lceil \frac{m+n}{2}\rceil}^{-1} \equiv 0 \text{ mod }  2^{\lfloor \frac{m+n}{2}\rfloor +t} \ee
and zero otherwise. Since $m\geq n+2$ equation \eqref{char2} has no solution  (and hence $S=0$)  unless $c_1=2^{n+t}c_3$ with $2\nmid   c_3$, in which  case \eqref{char2} becomes 
\be \label{char3}  \left( c_3A+     c_2 A_1 k_1 R_{\lceil \frac{m-n}{2}\rceil }R_{\lceil \frac{m+n}{2}\rceil}^{-1}    \right) 5^{vk}  \equiv -c_3B \text{ mod  } 2^{\lfloor \frac{m-n}{2}\rfloor }.    \ee
If no $v$ satisfies \eqref{char2} then plainly $S=0$. So assume that  \eqref{char2} has a solution   $v=v_0$ and take  $I_1=[v_0, v_0 +  2^{\lceil \frac{m-n}{2}\rceil -t -2})$.
Now  any other $v$ solving \eqref{char3}  must have
$$ 5^{vk}\equiv 5^{v_0k} \text{ mod } 2^{\lfloor \frac{m-n}{2}\rfloor }  \Rightarrow vk\equiv v_0k \text{ mod } 2^{\lfloor \frac{m-n}{2}\rfloor -2}  \Rightarrow v\equiv v_0 \text{ mod }  2^{\lfloor \frac{m-n}{2}\rfloor -t-2}. $$
So if  $m-n$ is even,  $I_1$
contains only the solution $v_0$ and
\be \label{SAeval}  S(A) = 2^{\lfloor \frac{m+n}{2}\rfloor +t}   \chi_1(5^{v_0}) \chi_2 (A5^{{v_0}k}+B). \ee
Observe that a solution $x_0=5^{v_0}$ or $x_0=-5^{v_0}$ of \eqref{defchar} corresponds to a solution $v_0$ to  \eqref{char2} when $k$ is even  and a solution $v_0$  to \eqref{char2} for $A$ or $-A$ respectively (both can not have solutions)  
if $k$ is odd.  The evaluation for $S$ follows at once  from \eqref{SAeval} and \eqref{trans1} or  \eqref{trans2}.
When $m-n$ is odd, $I_1$ contains two solutions $v_0$ and $v_0+ 2^{\lfloor \frac{m-n}{2}\rfloor -t-2}$  and
\begin{align*}  S(A) & =  2^{\lfloor \frac{m+n}{2}\rfloor +t}\chi_1(5^{v_0}) \left( \chi_2(h(5^{v_0})) +   \chi_1( 5^{ 2^{\lfloor \frac{m-n}{2}\rfloor -t-2}      }) 
\chi_2( A 5^{v_0k}  5^{k2^{\lfloor \frac{m-n}{2}\rfloor -t-2}     } +B)\right)  \\
 & = 2^{\lfloor \frac{m+n}{2}\rfloor +t} \chi_1(5^{v_0})\chi_2(h(5^{v_0})) \left( 1 + \chi_1(5^{2^{\lfloor \frac{m-n}{2}\rfloor -t-2}}) \chi_2(\xi ) \right) 
\end{align*}
where, since $3 \lfloor \frac{m-n}{2}\rfloor +n \geq m$ for $m\geq n+3$, 
\begin{align*}
\xi & =  A5^{v_0k}\left( 5^{k_12^{\lfloor \frac{m-n}{2}\rfloor -2}}-1\right) h(5^{v_0})^{-1} +1   \\ 
 & = A5^{v_0k}h(5^{v_0})^{-1} \left(  (1+R_{\lfloor \frac{m-n}{2}\rfloor}2^{\lfloor \frac{m-n}{2}\rfloor})^{k_1}-1\right) +1 \\
 & \equiv A5^{v_0k}h(5^{v_0})^{-1} \left( k_1R_{\lfloor \frac{m-n}{2}\rfloor}2^{\lfloor \frac{m-n}{2}\rfloor} + \binom{k_1}{2}R_{\lfloor \frac{m-n}{2}\rfloor}^2 2^{m-n-1} \right) +1 \text{ mod } 2^m \\
 & \equiv \left(  A_15^{v_0k}h(5^{v_0})^{-1} k_1R_{\lfloor \frac{m-n}{2}\rfloor} R_{\lfloor \frac{m+n}{2}\rfloor}^{-1} + \frac{1}{2}(k_1-1) 2^{\lfloor \frac{m-n}{2}\rfloor}    \right)R_{\lfloor \frac{m+n}{2}\rfloor} 2^{\lfloor \frac{m+n}{2}\rfloor}  +1 \text{ mod } 2^m \\
 & \equiv  \left( 1+ R_{\lfloor \frac{m+n}{2}\rfloor}  2^{\lfloor \frac{m+n}{2}\rfloor}  \right)^{A_15^{v_0k}h(5^{v_0})^{-1} k_1R_{\lfloor \frac{m-n}{2}\rfloor} R_{\lfloor \frac{m+n}{2}\rfloor}^{-1} + \frac{1}{2}(k_1-1)2^{\lfloor \frac{m-n}{2}\rfloor} } \text{ mod } 2^m  \\
 & = 5^{ \left(A_15^{v_0k}h(5^{v_0})^{-1} k_1R_{\lfloor \frac{m-n}{2}\rfloor}R_{\lfloor \frac{m+n}{2}\rfloor}^{-1}  + \frac{1}{2}(k_1-1)2^{\lfloor \frac{m-n}{2}\rfloor} \right) 2^{\lfloor \frac{m+n}{2}\rfloor -2}}.
 \end{align*}
Hence, setting 
$$    c_3 +   c_2    A_15^{v_0k}h(5^{v_0})^{-1} k_1R_{\lceil \frac{m-n}{2}\rceil} R_{\lceil \frac{m+n}{2}\rceil}^{-1} = \lambda 2^{ \lfloor \frac{m-n}{2}\rfloor}   $$
(only the parity of $\lambda$ will be used) and recalling that $c_2$ is odd, we have 
\begin{align*}
\chi_1(5^{2^{\lfloor \frac{m-n}{2}\rfloor -t-2}}) \chi_2(\xi ) & =e_{2^{\lceil \frac{m-n}{2}\rceil}} \left(c_3 +   c_2    A_15^{v_0k}h(5^{v_0})^{-1} k_1R_{\lfloor \frac{m-n}{2}\rfloor} R_{\lfloor \frac{m+n}{2}\rfloor}^{-1}  \right)    (-1)^{\frac{1}{2}(k_1-1)c_2}\\
 & = e_{2^{\lceil \frac{m-n}{2}\rceil}}\left(  c_2    A_15^{v_0k}h(5^{v_0})^{-1} k_1\left(R_{\lfloor \frac{m-n}{2}\rfloor} R_{\lfloor \frac{m+n}{2}\rfloor}^{-1} -   R_{\lceil \frac{m-n}{2}\rceil} R_{\lceil \frac{m+n}{2}\rceil}^{-1}    \right) \right)   (-1)^{\frac{1}{2}(k_1-1)+\lambda }.
\end{align*}
Since $1+R_{i+1}2^{i+1}=(1+R_i2^i)^2$ we have
$$ R_{i+1} = R_i + 2^{i-1}R_i^2 \equiv R_i+ 2^{i-1} \text{ mod } 2^{i+2}, $$
giving  $ R_i\equiv 3 \text{ mod } 4$ for $ i\geq 3, $ and
\begin{align*}  &  R_{\lfloor \frac{m-n}{2}\rfloor} R_{\lfloor \frac{m+n}{2}\rfloor}^{-1} -   R_{\lceil \frac{m-n}{2}\rceil} R_{\lceil \frac{m+n}{2}\rceil}^{-1}    \\
 \equiv  &   R_{\lfloor \frac{m+n}{2}\rfloor}^{-1}  R_{\lceil \frac{m+n}{2}\rceil}^{-1}  \left(   (R_{\lceil \frac{m-n}{2}\rceil}-2^{  \lceil \frac{m-n}{2}\rceil -2}) R_{ \lceil \frac{m+n}{2}\rceil} -   R_{\lceil \frac{m-n}{2}\rceil}(R_{  \lceil \frac{m+n}{2}\rceil} -    2^{  \lceil \frac{m-n}{2}\rceil +n-2} )   \right)  \text{ mod } 2^{ \lceil \frac{m-n}{2}\rceil } \\
\equiv &  (1-2^n) 2^{ \lceil \frac{m-n}{2}\rceil -2} \text{ mod }  2^{ \lceil \frac{m-n}{2}\rceil }.\end{align*}
From \eqref{char2} we have $ c_2    A_15^{v_0k}h(5^{v_0})^{-1} k_1\equiv -c_3 $ mod 4  and
$$ S(A)= 2^{\lfloor \frac{m+n}{2}\rfloor +t} \chi_1(5^{v_0})\chi_2(h(5^{v_0})) \left(1+ i^{(2^n-1)c_3} (-1)^{\frac{1}{2}(k_1-1)+\lambda}\right). $$
The result follows on writing 
$$ \frac{1+i^h}{\sqrt{2}} = \omega^h \left(\frac{2}{h}\right). $$

\vspace{2ex}
\noindent
{\bf  Small $m$ values: $t+2\leq m-n\leq 2t+4$}

Since $n+2(t+2)\geq m$ we have
\begin{align*}
A5^{\gamma  k} + B & = A_12^n (1+R_{t+2}2^{t+2} )^{\gamma k_1} +B \\
 & \equiv  (A+B)\left(1 + \gamma k_1 A_1 R_{t+2} (A+B)^{-1} 2^{t+n+2}\right) \text{ mod } 2^m\\ 
 & \equiv (A+B) \left(1+R_{t+n+2}2^{t+n+2}\right )^{\gamma k_1 A_1 (A+B)^{-1} R_{t+2}R_{t+n+2}^{-1}} \text{ mod } 2^m \\
 & =   (A+B) 5^{\gamma Ak (A+B)^{-1} R_{t+2}R_{t+n+2}^{-1}}.
\end{align*} 
Hence  $\chi_1(5^{\gamma}) \chi_2(A5^{\gamma k}+B) $ equals 
$$ \chi_2(A+B) e_{2^{m-2}}\left(\gamma \left(  c_1 (A+B) 
+ c_2Ak R_{t+2}R_{t+n+2}^{-1} \right) (A+B)^{-1}\right)$$
and $S(A) = 2^{m-2}\chi_2(A+B)$ if $C(1) \equiv 0$ mod $2^{m-2}$ and 0 otherwise.
Since $m-n\geq t+2$  the congruence $C(1)\equiv 0$ mod $2^{m-2}$ implies  $c_1=2^{t+n}c_3$
(with $c_3$ odd if $m-n>t+2$) and becomes
\be \label{keycong}  c_3(A+B) +c_2 A_1k_1  R_{t+2}R_{t+n+2}^{-1} \equiv 0 \text{ mod } 2^{m-n-t-2}. \ee
For $m-n=t+2$ or $t+3$ this will automatically hold (for both  $A$ and $-A$ when $k$ is odd)
and $S=2^{m-1}\chi_2(A+B)$ for $k$ even  and $\chi_1(-1)=1$,  and
$$S=2^{m-2}(\chi_2(A+B)+\chi_1(-1)\chi_2(-A+B))$$
for $k$ odd.
Further for $k$ odd and $m-n=2$ we have
 $-A+B\equiv  (1+2^{m-1})(A+B)$ mod $2^m$ with $\chi_2(1+2^{m-1})=-1$ and $S=2^{m-2}\chi_2(A+B)(1-\chi_1(-1))=2^{m-1}\chi_2(A+B)$ if $\chi_1(-1)=-1$
and zero otherwise.
Note when $m-n=t+2$ we have $c_1=2^{m-2}$ 
and $\chi_1(5)=1$ and when $m-n=t+3$  we have $c_1=2^{m-2}$ or $2^{m-3}$ and $\chi_1(5)=\pm 1$. 

Since $c_3B$ is odd
\eqref{keycong}  can not hold for both $A$ and $-A$ for $m-n>t+3$  and at most one of $S(A)$ or $S(-A)$ is non-zero. When $k$ is odd the congruence condition for $-A$ becomes $C(-1)\equiv 0$ mod $2^{m-2}$.

\end{document}